\documentclass[reqno,11pt]{amsart}
\usepackage{amsmath, amsfonts, amsthm, amssymb, setspace, textcomp, bbm, multirow, hyperref, mathtools, enumitem}
\usepackage{geometry}
\numberwithin{equation}{section}
\newtheorem{theorem}{Theorem}[section]
\newtheorem{lemma}{Lemma}[section]
\newtheorem{corollary}{Corollary}[section]

\newtheorem{remark}{Remark}[section]

\title{When Fourth Moments Are Enough}

\author{CHRIS JENNINGS-SHAFFER}
\thanks{Mathematical Institute, University of Cologne,
50931 K{\"o}ln, Germany.
jennichr@math.oregonstate.edu}

\author{Dane R. Skinner}
\thanks{Department of Mathematics, Oregon State University,
Corvallis, OR 97331.
skinner@onid.oregonstate.edu}

\author{Edward C. Waymire}
\thanks{Department of Mathematics, Oregon State University,
Corvallis, OR 97331.
waymire@math.oregonstate.edu}

\begin{document}

\begin{abstract}
This note concerns a somewhat innocent
question motivated by an observation
 concerning the use of Chebyshev 
bounds on sample estimates of $p$ in the binomial distribution
with parameters $n,p$.  Namely, what moment order produces
the best Chebyshev estimate of $p$?
If $S_n(p)$ has a binomial distribution with parameters $n,p$,
there it is readily observed that ${\rm argmax}_{0\le p\le 1}{\mathbb E}S_n^2(p) 
= {\rm argmax}_{0\le p\le 1}np(1-p) = \frac12,$ and
${\mathbb E}S_n^2(\frac12)  = \frac{n}{4}$.  
Rabi Bhattacharya observed that while the second moment Chebyshev sample size 
for a $95\%$ confidence estimate within $\pm 5$ percentage points is $n = 2000$, 
the fourth moment yields the substantially reduced polling requirement of
$n = 775$. Why stop at fourth moment? Is the argmax achieved at $p = \frac12$ for
higher order moments and, if so, does it help, and compute 
$\mathbb{E}S_n^{2m}(\frac12)$?   As captured by the title of this note, 
answers to these questions lead to a simple rule of thumb for best choice of 
moments in terms of an effective sample size for Chebyshev concentration inequalities.
\end{abstract}
\maketitle

\section{Introduction} 

This note concerns a somewhat innocent question motivated by an observation 
concerning the use of Chebyshev bounds on sample estimates of $p$ in the 
binomial distribution with parameters $n,p$. Namely, what moment order 
produces the best Chebyshev estimate of $p$? Chebyshev is arguably the most 
basic concentration inequality to occur in risk probability estimates and the 
use of second moments is a textbook example in elementary probability and 
statistics. Consider i.i.d. Bernoulli $0-1$ random variables $X_1,X_2,\dots,X_n$ 
with parameter $p\in[0,1]$, and let $S_n(p) = \sum_{j=1}^n(X_j-p)$. There it 
is readily observed that 
${\rm argmax}_{0\le p\le 1}{\mathbb E}S_n^2(p) = {\rm argmax}_{0\le p\le 1}np(1-p) = \frac12$.
It is also a well-known probability exercise to check that 4-th moment 
Chebyshev bounds improve the rate of convergence that can more generally be 
used for a proof of the strong law of large numbers, e.g. see (Bhattacharya and
Waymire, 2016; p.100).  Somewhat relatedly, Rabi Bhattacharya 
(personal communication) recently noticed, after a mildly tedious calculation 
to check ${\rm argmax}_{0\le p\le 1}{\mathbb E}S_n^4(p) = \frac12$, that the 
second moment Chebyshev bound is rather significantly improved by consideration 
of fourth moments as well.  In particular, while the second moment Chebyshev 
sample size for a $95\%$ confidence estimate within $\pm 5$ percentage points 
is $n = 2000$, the fourth moment yields the substantially reduced polling 
requirement of $n = 775$.  While the Chebyshev inequality is one among several 
inequalities used to obtain sample estimates, it is no doubt the simplest; see 
(Bhattacharya and Waymire, 2016) for comparison of fourth order Chebyshev to 
other concentration inequality bounds, and (Skinner, 2017) for numerical 
comparisons to higher order Chebyshev bounds.  
 
So why stop at fourth moments?  Is 
${\rm argmax}_{0\le p\le 1}\mathbb{E}S_n^{2m}(p) = \frac12$ for all $m,n$ and, 
if so, does it improve the estimate? Somewhat surprisingly we were not able to 
find a resolution of such basic questions in the published literature. In any 
case, with the argmax question resolved in part $(a)$ of the theorem below,
part $(b)$ provides a direct computation of $\mathbb{E}S_n^{2m}(\frac12)$. Part 
$(c)$ then provides a more readily computable version.   
\begin{theorem}\label{mainthm}
\begin{enumerate}[label=$(\alph*)$]
\item For all $m\ge 1$ and all $n$ sufficiently large,
${\rm argmax}_{0\le p\le 1}\mathbb{E}S^{2m}_n(p) = \frac12$. 
\item For all positive $m$ and $n$,
$\mathbb{E}S^{2m}_n(\frac12) = 
4^{-m}\sum_{\mu\in\pi(m), |\mu|\le m\wedge n}\binom{2m}{2\mu_1,\dots,2\mu_{|\mu|}}
\binom{n}{|\mu|}$, 
\item For all positive $m$ and $n$, 
${\mathbb E}S^{2m}_n(\frac12) = 2^{-2m-n}\sum_{k=0}^n \binom{n}{k} (2k-n)^{2m}$.
\end{enumerate}\noindent
Here $\pi(m)$ is the set of ordered integer partitions of $m$, also referred to 
as integer compositions, and $|\mu|$ denotes the number of parts of 
$\mu\in\pi(m)$. We refer to $|\mu|$ as the {\it size} of the partition $\mu$.
\end{theorem}

The equivalent calculus challenge is to show for fixed $m$ that for all 
sufficiently large $n$,
\begin{equation}\label{eqMoment}
\text{argmax}_{0\le p\le 1}\frac{d^{2m}}{dt^{2m}}(pe^{qt} + qe^{-pt})^n|_{t=0}
= \frac12.
\end{equation}
The example below illustrates the challenge to locating absolute maxima
for such polynomials (in $p$), especially to proofs by mathematical 
induction.  The proof given here is based on explicit combinatorial computation 
of ${\mathbb E}S_n^{2m}(p)$ in terms of ordered partitions of $2m$,
after introducing a few preliminary
lemmas.  The lemmas are relatively simple to check using the
statistical independence and identical distributions of the terms $X_i-p$ and $X_j-p$,
$i\neq j$, and make good exercises in calculus, probability, and number theory. 
However let us first observe that part $(a)$ of the theorem does not hold for $m > n$. 

\

\noindent {\bf Counter example to Theorem \ref{mainthm}$(a)$ for
(small)  $n< m$:}
Observe for $n = 1$ and $m=2$, the function
$${\mathbb E}S_1^4(p) = p - 4p^2 + 6p^3 - 3p^4, \quad 0\le p\le 1,$$
has a {\it minimum} at $p=\frac12$, with two local maxima
at $\frac12 \pm \frac{\sqrt{2}}{4}$.  In particular,
$${\rm argmax}_{0\le p\le 1}{\mathbb E}S_1^4(p) 
= \frac12 \pm \frac{\sqrt{2}}{4}.$$  
In particular, the polynomial is generally {\it not} unimodal.
So the restriction to sufficiently large $n$
is necessary for part $(a)$ of Theorem \ref{mainthm}. There is also the question of how large is sufficiently
large. We do not address this here, but computations suggest a bound along the
lines of $m\le c\cdot n^{\varepsilon}$, with $\varepsilon$ a little less than $\frac12$.
We let $m_n$ denote the largest value of $m$, dependent on $n$, such that
Theorem \ref{mainthm}$(a)$ holds for all $m\le m_n$. We leave it as an open problem
to determine an exact formula for $m_n$ and to determine a formula for 
${\rm argmax}_{0\le p\le1}\mathbb{E}_n^{2m}(p)$ for $m>m_n$.

\section{Proofs and Remarks}

Let $\pi(2m)$ denote the set of ordered partitions of $2m$. We we will use 
$|\mu| = k$ to denote the number of parts of $\mu$.  Finally, for 
$\mu\in\pi(2m)$, let 
\begin{equation*}
f_i(\mu,p) = pq^{\mu_i} + q(-p)^{\mu_i}, \quad 
0\le p\le 1, q=(1-p), 1\le i\le |\mu|.
\end{equation*}

\begin{lemma}
\label{basiclemma} Let $0\le p\le 1$ and $q = 1-p$. The following hold,
\begin{enumerate}[label=$(\alph*)$]
\item $S_n(p) = -^{\text{dist}} S_n(q)$,
\item $\mathbb{E}S^{2m}_n(p) = \mathbb{E}S^{2m}_n(q)$,
\item  $\mathbb{E}S^{2m}_n(p)
= \sum_{\mu\in\pi(2m)}\binom{n}{|\mu|} \binom{2m}{\mu_1,\dots,\mu_{|\mu|}}\prod_{i=1}^{|\mu|}f_i(\mu,p)$,
\item $\frac{d}{dp}\mathbb{E}S^{2m}_n(p)
= \sum_{\mu\in\pi(2m)}\binom{n}{|\mu|} \binom{2m}{\mu_1,\dots,\mu_{|\mu|}}\sum_{i=1}^{|\mu|}
f_i^\prime(\mu,p)\prod_{j\ne i}^{|\mu|}f_j(\mu,p)$.
\end{enumerate}
\end{lemma}

\begin{lemma}
\label{derivlem}
Let $\mu\in\pi(2m)$ and $1\le i\le |\mu|$. Then,
$$\frac{d}{dp}f_i(\mu,p) = 
q^{\mu_i}\left(1-\frac{p}{q}\mu_i\right) + (-1)^{\mu_i+1}
p^{\mu_i}\left(1-\frac{q}{p}\mu_i\right).$$
\end{lemma}

It now follows easily that 
\begin{equation}
\label{fi}
f_i\left(\mu,\frac12 \right) = 
\begin{cases}
2^{-\mu_i} & \text{for even}\  \mu_i,\\
0 &  \text{for odd}\  \mu_i,
\end{cases}
\end{equation}
\begin{equation}
\label{fiprime}
f^\prime_i\left(\mu,\frac12 \right) = 
\begin{cases}
0 & \text{for even}\  \mu_i,\\
-2(\mu_i-1)2^{-\mu_i} &  \text{for odd}\  \mu_i.
\end{cases}
\end{equation}

The keys to the following proof of Theorem \ref{mainthm} reside in 
(1) the parity conflicts between
\eqref{fi} and \eqref{fiprime} and (2) the expansion 
$(d)$ in Lemma \ref{basiclemma}, viewed as a polynomial in $n$.

\begin{proof}[Proof (of theorem)]
That $p=\frac12$ is a critical point follows from $(d)$ of Lemma \ref{basiclemma}
together with \eqref{fi} and \eqref{fiprime} by examining the terms
$f_i^\prime(\mu,\frac12)\prod_{j\ne i}^{|\mu|}f_j(\mu,\frac12)$.
In particular, for partitions of $2m$, if $\mu_i$ is odd then there must be a 
$j\neq i$ such that $\mu_j$ is odd as well. To see that $p=\frac12$ is an 
absolute maximum, the trick is to observe that for $0\le p < \frac12 < q$, 
the leading coefficient of $\frac{d}{dp}\mathbb{E}S_n^{2m}(p)$, viewed as a 
polynomial in $n$, is obtained at the $m$-part composition $\mu = (2,2,\dots,2)$ 
of $2m$. Namely, it is obtained from
$\binom{n}{m}\binom{2m}{2,2,\dots,2}m(q^2-p^2)(pq)^{m-1}$, 
and therefore is positive for all $p<\frac12$. Thus, for sufficiently large $n$, 
$$\frac{d}{dp}\mathbb{E}S_n^{2m}(p) > 0, \quad\mbox{for } 0\le p < 1/2.$$
In view of the symmetry expressed in
$(b)$ of Lemma \ref{basiclemma}, it follows that $p =\frac12$ is the unique
global maximum.

For part $(b)$ of the theorem one simply computes from independence,
writing $\tilde{X}_i = X_i-\frac12, i =1,2,\dots, n$. In particular,
$\tilde{X}_i = \pm\frac12$ with equal probabilities. So,
for $m\ge 1$,
\begin{align*}
\mathbb{E}S_n^{2m}\left(\frac12\right) 
&= 
\sum_{1\le j_1,\dots,j_{2m}\le n}{\mathbb E}\prod_{i=1}^{2m}\tilde{X}_{j_i}
\\
&=
\sum_{2m_1+\cdots +2m_n = 2m}\prod_{i=1}^{n}{\mathbb E}\tilde{X}_{i}^{2m_i}
\\
&= 
\sum_{k=1}^{m\wedge n}\sum_{2m_1+\cdots +2m_n = 2m, \#\{j:m_j\ge 1\}=k}\prod_{i=1}^{n}4^{-m_i}
\\
&= 
\sum_{k=1}^{m\wedge n}\binom{n}{k}\sum_{\mu=(\mu_1,\dots,\mu_k)\in\pi(m)}
\binom{2m}{2\mu_1,\dots,2\mu_k}4^{-m}.
\end{align*}
Here one adopts the convention that a sum over an empty set is
zero so that if there are no partitions $\mu$ of $m$ with $|\mu|=k$
then the indicated sum is zero for this choice of $k$. So 
nonzero contributions to the  sum are provided by ordered partitions $\mu$
of size $|\mu|\le m\wedge n$.

To simplify the computation in terms of ordered partitions $(b)$ one
may proceed as follows to obtain the formula in $(c)$. We instead compute
$\mathbb{E}S_n^{2m}(\frac12)$ as the $2m$-th moment of $S_n(\frac12)$ as
given in \eqref{eqMoment}. By the binomial theorem, we have that
\begin{align*}
\mathbb{E}S_n^{2m}\left(\frac12\right)
&=
\frac{d^{2m}}{dt^{2m}}
\left[\left( 
	\frac{e^{\frac{t}{2}}}{2} 
	+
	\frac{e^{-\frac{t}{2}}}{2}
\right)^n\right]_{t=0}
=
\frac{d^{2m}}{dt^{2m}}
\left[
2^{-n}
\sum_{k=0}^n \binom{n}{k}
e^{\frac{t}{2}(2k-n)}
\right]_{t=0}
\\&=
2^{-n-2m}\sum_{k=0}^n \binom{n}{k} (2k-n)^{2m}.
\end{align*}
\end{proof}

\begin{remark}
A linear recurrence in $m$ is possible 
to aid the pre-asymptotic (in $n$) computation
of $\mathbb{E}S_n^{2m}(\frac12)$.
Namely, 
\begin{align}
\label{recur}
\mathbb{E}S_n^{2m+2\ell+2}\left(\frac12\right)
&=
\sum_{j=0}^\ell
c_j 2^{2j-2\ell-2}\mathbb{E}S_n^{2m+2j}\left(\frac12\right),
\end{align}
where
$\ell=\left\lfloor\frac{n-1}{2}\right\rfloor$,
$a_k = (2k-n)^2$, and $(c_0,c_1,\dotsc,c_\ell)$ is the unique solution to
\begin{align*}
\begin{pmatrix}
a_0^0&a_0^1&\dots&a_0^\ell\\
a_1^0&a_1^1&\dots&a_1^\ell\\
\vdots
\\
a_\ell^0&a_\ell^1&\dots&a_\ell^\ell\\
\end{pmatrix}
\begin{pmatrix}
c_0\\c_1\\ \vdots\\ c_\ell
\end{pmatrix}
&=
\begin{pmatrix}
a_0^{\ell+1}\\ a_1^{\ell+1}\\ \vdots\\ a_\ell^{\ell+1}
\end{pmatrix}.
\end{align*}
To see this, write
\begin{align*}
\mathbb{E}S_n^{2m}\left(\frac12\right)
&=
	2^{-2m-n+1}\sum_{k=0}^{\ell} 
	\binom{n}{k} (2k-n)^{2m}	.
\end{align*}
Then \eqref{recur} follows  since
\begin{align*}
&\mathbb{E}S_n^{2m+2\ell+2}\left(\frac12\right)
-
\sum_{j=0}^\ell
c_j 2^{2j-2\ell-2}\mathbb{E}S_n^{2m+2j}\left(\frac12\right)
\\
&=
	2^{-2m-2\ell-n-1}\sum_{k=0}^{\ell}\binom{n}{k}a_k^{m+\ell+1}
	-
	\sum_{j=0}^{\ell}c_j 2^{-2m-2\ell-n-1}
	\sum_{k=0}^{\ell}\binom{n}{k}a_k^{m+j}
\\
&=
	2^{-2m-2\ell-n-1}\sum_{k=0}^{\ell}\binom{n}{k}a_k^{m}
	\left(a_k^{\ell+1}-\sum_{j=0}^{\ell}c_ja_k^j \right)	
=0.
\end{align*}

\end{remark}

\

\

For the application to statistical estimation one may 
combine Theorem \ref{mainthm}
 with  Chebyshev's inequality to obtain,
\begin{corollary} For $\epsilon > 0$, we have that
$
P\left(|\frac{1}{n}S_n(p)| > \epsilon\right) \le \min_{1\le m\le m_n}
\left(
\frac{\sqrt[2m]{{\mathbb E}S^{2m}_n(\frac12)}}
{n\epsilon}\right)^{2m}.
$
\end{corollary}
Noting the scaling invariance
$
\text{argmax}_{0\le p \le 1}\mathbb{E}S_n^{2m}(p) = 
\text{argmax}_{0\le p\le 1}\mathbb{E}\frac{S_n^{2m}(p)}{n^m},
$
and 
${\mathbb E}Z^{2m} = 2^{-m}\frac{(2m)!}{m!}$ for the standard normal random 
variable $Z$, in the limit ``$n \to\infty, \epsilon \to 0, n\epsilon^2\to \tilde{n}$'' 
one has
\begin{equation*}
B_m := 
\mathbb{E}\frac{S^{2m}_n(\frac12)}{n^{2m}\epsilon^{2m}}
= 
\mathbb{E}
\frac{\left(\frac{S_n(\frac12)}{\sqrt{n/4}}\right)^{2m}}{n^{2m}\epsilon^{2m}}
\left(\frac{n}{4}\right)^m 
\to 
2^{-2m}\tilde{n}^{-m}{\mathbb E}Z^{2m} = 2^{-3m}\frac{(2m)!}{m!}\tilde{n}^{-m}.
\end{equation*}
In particular, one may ask for the best choice of $m$ for large $n$, i.e,
in the above limit as $n\to\infty,\epsilon\downarrow 0, n\epsilon^2\to \tilde{n}$. 
The quantity $\tilde{n} = n\epsilon^2$ denotes
an {\it effective sample size} in the sense of the risk assessment defined
by $P(|S_n(p)| > n\epsilon) < \epsilon$; see (Duchi et al. 2013) for
an introduction of this artful terminology in a much broader context.
Observe that in the limit of large $n$
\begin{equation*}
\lim_{n\to\infty, \epsilon\downarrow 0, n\epsilon^2 = \tilde{n}}
\frac{B_{m+1}}{B_m} = \frac{2m+1}{4\tilde{n}}
\begin{cases}
\le 1\\
= 1\\
\ge 1
\end{cases}
\end{equation*}
if and only if 
\begin{equation*}
m 
\begin{cases}
\le 2\tilde{n} -\frac12\\
= 2\tilde{n} -\frac12\\
\ge 2\tilde{n} -\frac12.
\end{cases}
\end{equation*}

\

\

The take-away is perhaps best summarized in terms of the following
informally interpreted optimal estimation principle.

\

\

\noindent {\bf Approximate Rule of Thumb:} 
{\it For large $n$ the optimal moment order $2m$
for the Chebyshev bound is quadruple the effective sample size.
In particular,  the fourth moment is optimal for a one unit
effective sample size!}


\section{References}

\noindent Bhattacharya, R.N., and  Waymire, E.C. (2016),
``A Basic Course in Probability Theory'', 2nd ed.,
Universitext,
Springer, NY.

\ 

\noindent Duchi, J. and Wainwright, M. J. and Jordan, M. I. (2013),
``Local privacy and minimax bounds: Sharp rates for probability estimation'', 
in Advances in Neural Information Processing Systems,
1529-1537.

\

\noindent Skinner, Dane, 2017, ``Concentration of Measure Inequalities'', 
Master of Science, Oregon State University
 \ 



 \end{document}